\documentclass[12pt]{article}
\usepackage{amsmath,amssymb,amstext,amsthm,bm,amsfonts}
\usepackage{url}
\usepackage{color,graphicx,graphics,xcolor}
\usepackage{hyperref}
\usepackage{mathrsfs}
\usepackage{listings}
\usepackage{appendix}

\usepackage{upref,amsxtra,exscale}
\usepackage{cite}
\usepackage{hyperref}

\usepackage {comment}

\usepackage{chngcntr}
\counterwithin*{equation}{section}

\newtheorem{thm}{Theorem}[section]
\newtheorem{defi}[thm]{Definition}

\newtheorem{cor}[thm]{Corollary}
\newtheorem{lem}[thm]{Lemma}
\newtheorem{rem}[thm]{Remark}

\newtheorem{example}[thm]{Example}

\newcommand{\R}{\mathbb{R}}

\def\qq#1{\qquad \mbox{#1}\quad}
\newcommand{\D}{\displaystyle}

\newcommand{\al}{\alpha}
\newcommand{\be}{\beta}
\newcommand{\De}{\Delta}
\newcommand{\de}{\delta}

\newcommand{\e}{\varepsilon}

\newcommand{\ga}{\gamma}

\newcommand{\la}{\lambda}

\newcommand{\na}{\nabla}

\newcommand{\Om}{\Omega}
\newcommand{\Omb}{\overline{\Om}}
\newcommand{\p}{\partial}
\newcommand{\s}{\sigma}
\newcommand{\te}{\theta}

\newcommand{\vf}{\varphi}

\title{$L^\infty$ a-priori estimates for subcritical semilinear elliptic equations with  a Carathéodory nonlinearity}
\author{Rosa Pardo\\
Departamento de An\'alisis Matem\'atico y Matem\'atica Aplicada\\ Universidad Complutense de Madrid, 28040 
Madrid,
Spain\\
rpardo@ucm.es}

\date{}

\begin{document}
\maketitle
\begin{abstract}
We present  new $L^\infty$ {\it a priori} estimates for weak solutions of a wide class of subcritical  elliptic equations in bounded domains.   No hypotheses on the sign of the solutions, neither of the non-linearities are required. This method is based in Gagliardo-Nirenberg and Caffarelli-Kohn-Nirenberg interpolation inequalities.

Let us consider a semilinear boundary value problem
$ -\Delta u= f(x,u),$ in $\Omega,$ with Dirichlet boundary conditions, where $\Omega \subset \R ^N $, with $N> 2,$ is a bounded smooth domain, and $f$  is  a subcritical Carathéodory non-linearity.	
We provide $L^\infty$ a priori estimates for weak  solutions, in terms of their $L^{2^*}$-norm, where $2^*=\frac{2N}{N-2}\ $ is the critical Sobolev exponent. 

By a subcritical non-linearity we mean, for instance,
$|f(x,s)|\le |x|^{-\mu}\, \tilde{f}(s),$
where $\mu\in(0,2),$  and $\tilde{f}(s)/|s|^{2_{\mu}^*-1}\to 0$ as $|s|\to \infty$,  here $2^*_{\mu}:=\frac{2(N-\mu)}{N-2}$ is the critical Sobolev-Hardy exponent.  Our non-linearities includes non-power non-linearities. 

In particular we prove that when $f(x,s)=|x|^{-\mu}\,\frac{|s|^{2^*_{\mu}-2}s}{\big[\log(e+|s|)\big]^\be}\,,$ with $\mu\in[1,2),$  
then, for any $\e>0$ there exists a constant $C_\e>0$ such that for any solution $u\in H^1_0(\Om)$, the following holds
$$
\Big[\log\big(e+\|u\|_{\infty}\big)\Big]^\be\le
C _\e \, \Big(1+\|u\|_{2^*}\Big)^{\, (2^*_{\mu}-2)(1+\e)}\, .
$$
\end{abstract}
\noindent\textbf{MSC2020:} {Primary 35B45; Secondary 
35J75  % Singular elliptic equations 
35J25, %(bvp ell eq), 
35J60, 35J61, 35B33,  35B65.
% 35J60   	Nonlinear elliptic equations
% 35J61   	Semilinear elliptic equations
}

\noindent\textbf{Keywords:} A priori estimates, subcritical non-linearities, $L^\infty$  a priori bounds, changing sign weights, singular elliptic equations

\maketitle

\section{Introduction}
Let us consider the following semilinear boundary value problem
\begin{equation}
\label{eq:ell:pb} 
-\Delta u= f(x,u), \quad\mbox{in }\Om, \qquad  u= 0,\quad \mbox{on }\p \Om,
\end{equation}
where $\Om \subset \R^N $,  $N> 2,$ is a bounded, connected, open subset with $C^{2}$ boundary $\p \Om$, 
and the non-linearity $f:\bar{\Om}\times\R\to \R$ is Carathéodory function (that is, the mapping $f(\cdot,s)$ is measurable for all $s\in\R$,  and the mapping $f(x,\cdot)$ is continuous for almost all $x\in\Om$),  
and {\it subcritical} (see definition \ref{def:sub}). 

\medskip

We analyze the effect of the smoothness of the non-linearity $f=f(x,u)$ on the $L^\infty(\Omega)$ a priori estimates of  {\it weak solutions} to \eqref{eq:ell:pb}. 
Degree theory combined with a priori bounds in the sup-norm of solutions of parametrized versions of \eqref{eq:ell:pb}, is a very classical topic in elliptic equations, posed 
by Leray and Schauder in\cite{Leray-Schauder}. It provides a great deal of information about existence of  solutions and the structure of the  solution set. This study is usually focused on positive classical solutions, see the classical references of de Figueiredo-Lions-Nussbaum, and of Gidas-Spruck \cite{deFigueiredo-Lions-Nussbaum, Gidas-Spruck}, see also  \cite{Castro_Pardo_RMC, Castro-Pardo_DCDS}.

\bigskip

A natural question concerning  the class  of  solutions is the following one, 
\begin{enumerate}
\item[(Q1)] {\it  those $L^\infty(\Omega)$ estimates can be extended to a bigger class of solutions, in particular to {\it weak solutions} (with possibly changing sign solutions)?. }
\end{enumerate}
Another natural question with respect to  the class of non-linearities, can be stated as follows,
\begin{enumerate}
\item[(Q2)] {\it  those estimates can be extended to a bigger class of non-linearities, in particular  to non-smooth non-linearities (with possibly changing sign weights)?. }
\end{enumerate}

\medskip

\noindent	In this paper we provide sufficient conditions guarantying uniform $L^\infty(\Omega)$ a priori estimates for any $u\in H_0^{1}(\Om)$ weak solution to \eqref{eq:ell:pb}, in terms of their $L^{2^*}(\Omega)$ bounds, in the class of Caratheodory subcritical generalized problems. 
In this class, we  state that any set of weak  solutions uniformly $L^{2^*} (\Om)$ a priori bounded is universally $L^\infty (\Om)$ a priori bounded. Our theorems allow changing sign weights, and singular weights, and also apply to changing sign solutions.

\bigskip

Problem \eqref{eq:ell:pb} with $f(x,s)=|x|^{-\mu} |s|^{p-1}s ,$  $\mu>0$, is known as Hardy's problem, due to its relation with the Hardy-Sobolev inequality.
The Caffarelli-Kohn-Nirenberg interpolation inequality for radial singular weights \cite{Caffarelli-Kohn-Nirenberg}, states that whenever $0\le\mu\le 2$, 
\begin{equation}\label{def:2*:mu}
2_{\mu}^*:=\frac{2(N-\mu)}{N-2},
\end{equation} 
is the {\it critical} exponent of the Hardy-Sobolev embedding $\qquad
H_0^1(\Om)\hookrightarrow\ $ $ L^{2_\mu^*}(\Om, |x|^{-\mu})\,
$. 
Using variational methods, one obtains the existence of a nontrivial solution to \eqref{eq:ell:pb} in $H_0^1(\Om)$  whenever $1 < p  <2_\mu^*-1$.
For the case $0<\mu< 2$, using a Pohozaev type identity, we have that for $p\ge 2_\mu^*-1$ there is no solution to Hardy's problem in star-shaped domains with respect to the origin. 
But, there exist positive solutions for the problem with $p=2_\mu^*-1$
depending on the geometry of the domain $\Om$, see \cite{Jannelli-Solimini} and \cite{Caldiroli-Malchiodi}.

\smallskip

If $\mu\ge 2,$ it is known that Hardy's problem has no positive solution in any domain $\Omega$ containing the origin, see \cite{Gidas-Spruck2}, \cite[Lemma 6.2]{Bidaut-Veron}, \cite{Dupaigne-Ponce}.

\bigskip

Usually the term subcritical non-linearity is reserved for power like non-linearities. Next, we expand this concept:

\begin{defi}
\label{def:sub}
By a {\it subcritical non-linearity} we mean that $f$ satisfies one of the following two growth conditions:

\begin{enumerate}				
\item[{\rm (H0) }] 
\begin{equation}\label{f:growth:rad}
|f(x,s)|\le |x|^{-\mu}\, \tilde{f}(s),
\end{equation}
\end{enumerate}
where $\mu\in(0,2),$  and $\tilde{f}:\R\to[0,+\infty)$   is continuous and satisfy  
\begin{equation}\label{f:sub:rad}
\tilde{f}(s)>0 \ \text{for}\ |s|>s_0,\ \text{ and }\  \lim_{|s| \to \infty} \ \frac{\tilde{f}(s)}{|s|^{2_{\mu}^*-1}}=0;
\footnote{
Observe that $2_{\mu}^*>2$ for $\mu\in(0,2)$. Let $a(x)=|x|^{-\mu}$,  $a\in L^p(\Om)$ for any $p<N/\mu$. Moreover, the sharp Caffarelli-Kohn-Nirenberg  inequality implies that if $u\in H_0^1(\Om)$, then
$f(\cdot,u)\in L^{\frac{2N}{N+2}}(\Om)$
}
\end{equation}
or
\begin{enumerate}				
\item[{\rm (H0)'}] 
\begin{equation}\label{f:growth}
|f(x,s)|\le |a(x)|\, \tilde{f}(s)
\end{equation}
\end{enumerate}
where $a\in L^r(\Om)$ with $r>N/2$, and  $\tilde{f}:\R\to [0,+\infty)$ is continuous and satisfy   
\begin{equation}\label{f:sub}
\tilde{f}(s)>0 \ \text{for}\ |s|>s_0,\ \text{ and }\  \lim_{s \to \pm\infty}  \frac{\tilde{f}(s)}{|s|^{2_{N/r}^*-1}}=0,
\end{equation}
where 
\begin{equation}\label{def:2:r}
2_{N/r}^*:=\frac{2^*}{r'}=2^*\left(1-\frac{1}{r}\right),
\end{equation}
and where $r'$ is the conjugate exponent of $r$,  $1/r+1/r'=1.$	
\footnote{Since $r>N/2,$ then  $2_{N/r}^*>2.$
Moreover, thanks to Sobolev embeddings, for any $u\in H_0^{1}(\Om)$,
\begin{align}		
% u\in L^{{2^*}}(\Om) & \qq{with} 2^*r'=\frac{2Nr}{(N-2)(r-1)},\\ 
\tilde{f}(u)\in L^{\frac{2^*}{2_{N/r}^*-1}}(\Om) & \qq{with} \frac{2_{N/r}^*-1}{2^*}=\frac12+\frac1N-\frac{1}{r}, \nonumber	
\\ 
\text{and}\ f(\cdot,u)\in L^{\frac{2N}{N+2}}(\Om). \nonumber	
\end{align}
}		
\end{defi}

\medskip

Our analysis shows that non-linearities satisfying either (H0): \eqref{f:growth:rad}-\eqref{f:sub:rad} (either (H0)': \eqref{f:growth}-\eqref{f:sub}),  widen the  class of subcritical non-linearities to non-power non-linearities, sharing with power like non-linearities properties such as $L^\infty$ a priori estimates. 
Our definition of a subcritical non-linearity includes  non-linearities such as 
$$
f^{(1)}(x,s):=
\frac{|x|^{-\mu}|s|^{2_{\mu}^*-2}s}{\big[\log(e+|s|)\big]^\al},
\ \ \text{or}\ \ f^{(2)}(x,s):= \frac{a(x)|s|^{2_{N/r}^*-2}s}{\Big[\log\big[e+\log(1+|s|)\big]\Big]^\al},
$$
for any $\al>0$, and  $\mu\in(0,2)$, or any $a\in L^r(\Om),$ with $r>N/2$.

\medskip

In particular, if $f(x,s)=f^{(1)}(x,s)$ with $\mu\in [1,2),$ then for any $\e>0$ there exists a constant $C>0$ depending  on $\e$,  $\mu$, $N$,  and $\Om,$ such that for any  $u\in H_0^{1}(\Om)$  solution   to \eqref{eq:ell:pb}, the following holds:
\begin{equation*}%\label{bdds:2}
\Big[\log\big(e+\|u\|_{\infty}\big)\Big]^\al\le
C \ 
\Big(1+\|u\|_{2^*}\Big)^{\, (2_{\mu}^*-2)(1+\e)},
\end{equation*}
and where  $C$  is independent of the  solution $u$, see Theorem \ref{th:rad}. Related results concerning those non-power non-linearities can be found in \cite{Damascelli-Pardo} for the $p$-laplacian case,   \cite{Clapp-Pardo-Pistoia-Saldana} analyzing what happen when $\al\to 0$, \cite{Mavinga-Pardo17_JMAA} for systems, \cite{Pardo-Sanjuan} for the radial case, and \cite{Pardo19a, Pardo19b} for a summary. 

\medskip

Moreover, if $f(x,s)=f^{(2)}(x,s)$ with $a\in L^r(\Om)$ for $r\in(N/2,N]$, then for any $\e>0$ there exists a constant $C>0$ depending only on $\e,\ \Om,$ $r$ and $N$ such that for any  $u\in H_0^{1}(\Om)$  solution   to \eqref{eq:ell:pb}, the following holds:
\begin{equation*}%\label{bdds:1}
\Big[\log\big[e+\log\big(1+\|u\|_{\infty}\big)\big]\Big]^\al\le
C \|a\|_r^{\ 1+\e} \ \Big(1+\|u\|_{2^*}\Big)^{\, (2_{N/r}^*-2)(1+\e)},
\end{equation*}
where  $C$  is independent of the  solution $u$, see Theorem \ref{th:apriori:cnys}.

\begin{defi}
\label{def:w:sol}
By a {\it solution} we mean a weak solution $u \in H_0^{1}(\Om)$  such that $f(\cdot,u)\in L^{\frac{2N}{N+2}}(\Om),$ and
\begin{equation}\label{weak:sol}
\int_{\Om}\nabla u\nabla \vf=\int_{\Om} f(x,u)\vf, \quad\quad \forall \vf\in H_0^1(\Om).
\end{equation}
\end{defi}

By an estimate of Brezis-Kato \cite{Brezis-Kato}, based on Moser's iteration technique \cite{Moser}, and according to elliptic regularity,   any weak solution to \eqref{eq:ell:pb} with a Caratheodory subcritical non-linearity is a continuous function, and in fact is an strong solution, see Lemma \ref{lem:reg} and Lemma \ref{lem:reg:rad}.

\bigskip

Joseph and Lundgren in \cite{Joseph-Lundgren} show us that those $L^\infty$ a priori estimates are not applicable for $L^1$- weak solutions, neither to super-critical non-linearities. 
\begin{defi}
\label{def:L1:w:sol}
We will say that a function $u$ is an  {\it $L^1$- weak solution} to \eqref{eq:ell:pb} if 	
$$
u\in L^1(\Om),\qquad f(\cdot,u)\, \de_\Om \in   L^1(\Om)
$$
where $\de_\Om (x):=dist(x,\p\Om)$ is the distance function with respect to the boundary, and 
\begin{equation*}%\label{eq:L1:sol}
\int_\Om \Big(u  \Delta \varphi +  f(x,u)\varphi\Big) \, dx=0, \qq{for all} \varphi \in C^2(\Omb ),\quad \varphi\big|_{\p\Om}=0. 
\end{equation*}
\end{defi}

They posed the study of singular solutions. Working on non-linearities such as $f(s):=e^s$ or $f(s):=(1+s)^p$, they consider the following BVP depending on a multiplicative parameter $\la\in\R$,
\begin{equation}\label{pde:la} 
-\Delta u=\la f(u), \quad\mbox{in }\Om, \qquad  u= 0,\quad \mbox{on }\p \Om,
\end{equation}
and look for classical radial positive solutions in  the unit ball $B_1$.
They obtain singular solutions as limit of classical solutions.

In particular, they obtain the explicit weak solution
\begin{equation*}%\label{sing:sol:exp}
u^{*}_1(x):=\log \frac{1}{|x|^2},
\quad u^{*}_1\in H_0^1(B_1),
\end{equation*}
to \eqref{BVP:sing:w}, when $N>2$, $\la=2(N-2)$,  and $f(s):=e^s$, see \cite[p. 262]{Joseph-Lundgren}. 

They also found the explicit $L^1${\textit -weak solution} 
\begin{equation*}%\label{sing:sol:pol}
u^{*}_2(x):= \left(\frac{1}{|x|}\right)^\frac{2}{p-1}-1, 
\ \ \text{with}\  p>\frac{N}{N-2},\ N>2,
\quad u^{*}_2\in W_0^{1,\frac{N}{N-1}}(B_1),
\end{equation*}
to \eqref{BVP:sing:w}, where $f(s):=(1+s)^p$, and $\la=\frac{2}{p-1} \big(N-\frac{2p}{p-1}\big)>0$,  
see \cite[(III.a)]{Joseph-Lundgren}. It holds that $u^{*}_2\in H_0^1(B_1)$ only when $p>2^*-1$.
So,  in the subcritical range   $u^{*}_2$ is a singular $L^1$-weak
solution, not in $H^1$.

\medskip

Let us focus on BVP with radial singular weights,
\begin{equation}\label{BVP:sing:w}
-\Delta u=\la |x|^{-\mu} (1+u)^p, \quad\mbox{in }\Om, \qquad  u= 0,\quad \mbox{on }\p \Om,
\end{equation}
with  $N>2$, $\mu<2$ and $p>1$. It can be checked that 
\begin{equation*}%\label{sing:sol:pol:rad}
u^{*}_3(x):= \left(\frac{1}{|x|}\right)^\frac{2-\mu}{p-1}-1,
\ \ \text{with}\ p>\frac{N-\mu}{N-2},
\qquad u^{*}_3\in W_0^{1,\frac{N}{N-1}}(B_1),
\end{equation*}
and $u^{*}_3$ is an $L^1${\textit -weak solution} to \eqref{BVP:sing:w},
for $\la=\frac{2-\mu}{p-1} \big(N-2-\frac{2-\mu}{p-1}\big)>0$. It also holds that  $u^{*}_3\in H_0^1(B_1)$ only when $p>2_\mu^*-1$.
So,  in the subcritical range  $u^{*}_3$ is a singular $L^1$-weak
solution  to \eqref{BVP:sing:w}, not in $H^1$.

\bigskip

Those examples of radially symmetric singular solutions  to BVP's on spherical domains,  solve either super-critical problems ($u_1^*$) or are  $L^1$-{\textit weak solutions} not in $H_0^1(\Om)$ ($u_2^*$ and $u_3^*$).
Consequently,  we restrict our study for $u\in H_0^{1}(\Om)$ {\it weak solutions} to \eqref{eq:ell:pb}, in the class of subcritical generalized problems. It is natural to ask for uniform $L^\infty$ a priori estimates   over non power non-linearities in non-spherical domains.

\bigskip

To state our main results, for a non-linearity $f$ satisfying \eqref{f:growth:rad}-\eqref{f:sub:rad}, let us define
\begin{equation}\label{rad:def:h}
h(s)=h_{\mu}(s):=\frac{|s|^{2_{\mu}^*-1}}{\displaystyle 
\max_{\{-s,s\}} \tilde{f}}, \qq{for} |s|> s_0.
\end{equation}
And for a non-linearity $f$ satisfying \eqref{f:growth}-\eqref{f:sub}, let us now define

\begin{equation}\label{def:h}
h(s)=h_{N/r}(s):=\frac{|s|^{2_{N/r}^*-1}}{\displaystyle 
\max\big\{\tilde{f}(-s),\tilde{f}(s)\big\}}\qq{for} |s|> s_0.
\end{equation}
By sub-criticallity, (see \eqref{f:sub:rad} or \eqref{f:sub} respectively), 
\begin{equation}\label{h:inf}
h(s)\to\infty \qq{as} s\to\infty.
\end{equation}

\medskip

\noindent Let $u$ be a solution to \eqref{eq:ell:pb}. We estimate $h\big(\|u\|_{\infty}\big)$, in terms of its $L^{2^*}$-norm. 
This result is robust, and holds for positive, negative and changing sign non-linearities, and also for positive, negative and changing sign solutions.

As an immediate consequence, as soon as we have a universal {\it a priori} $L^{2^*}$- norm for weak  solutions in $H_0^{1}(\Om),$ 
%(see footnote \ref{foot:reg}), 
then solutions are a priori universally bounded  in the $L^{\infty}$- norm,  see Corollary \ref{cor:apriori:cnys}.
\bigskip

\medskip

This paper is organized in the following way. 
In  Section \ref{sec:main:I}, using Gagliardo --Nirenberg inequality, we analyze the case when $a\in L^r(\Om)$ with $r>N/2$, see Theorem \ref{th:apriori:cnys}. In  Section \ref{sec:main:II}, we analyze the more involved case of a radial singular weight, see Theorem \ref{th:rad}. It needs the  Caffarelli-Kohn-Nirenberg inequality.

\section{Carathéodory non-linearities}
\label{sec:main:I}

In this section, assuming that $f$ satisfy the subcritical growth condition, we state our first main result concerning Carathéodory non-linearities, see  Theorem \ref{th:apriori:cnys}.  

\smallskip

We first collect a regularity Lemma for any weak solution to \eqref{eq:ell:pb} with a  non-linearity  of polynomial critical growth.

\begin{lem}[Improved regularity]
\label{lem:reg} 

Assume that $u\in H_0^1(\Om)$ weakly solves \eqref{eq:ell:pb} for a Carathéodory non-linearity $f:\bar{\Om}\times\R\to\R$ with polynomial critical growth
\begin{equation}
\label{f:growth:2}
|f(x,s)|\le |a(x)|(1+|s|^{2_{N/r}^*-1}),\quad\text{with}\ a\in L^r(\Om),\quad N/2<r\le\infty.
\end{equation}
Then, the following hold: 
\begin{enumerate}
\item[\rm (i)] If $r<N,$ then $u\in C^{\nu}(\overline{\Om})\cap W^{2,r}(\Om)$ for $\nu=2-\frac{N}{r}\in (0,1)$.
\item[\rm (ii)] If $r=N,$ then $u\in C^{\nu}(\overline{\Om})\cap W^{2,r}(\Om)$ for any $\nu<1$.
\item[\rm (iii)] If $N<r<\infty,$ then $u\in C^{1,\nu}(\overline{\Om})\cap W^{2,r}(\Om)$ for $\nu=1-\frac{N}{r}\in (0,1)$.
\item[\rm (iv)] If $r=+\infty,$ then $u\in C^{1,\nu}(\overline{\Om})\cap W^{2,p}(\Om)$ for any $\nu<1$ and any $p<\infty$.
\end{enumerate}  

\end{lem} 
\begin{proof}
Let $u\in H_0^1(\Om)$ be a solution to \eqref{eq:ell:pb}. Since an estimate of Brezis-Kato \cite{Brezis-Kato}, %based on Moser's iteration technique \cite{Moser}, 
if
\begin{equation}
\label{f:growth:3}
|f(x,u)|\le b(x)(1+|u|),\qq{with} 0\le b\in L^{N/2}(\Om),
\end{equation}
then,  $u\in L^q(\Om)$ for any $q<\infty$ (see \cite[Lemma B.3]{Struwe}).

Assume that $f$ satisfies \eqref{f:growth:2}, then assumption \eqref{f:growth:3} is satisfied with
$$
b(x)=\frac{|a(x)|(1+|u|^{2_{N/r}^*-1})}{1+|u|}\le C\, |a(x)|(1+|u|^{2_{N/r}^*-2})\in L^{N/2}(\Om).
$$
Consequently, $u\in L^q(\Om)$ for any $q<\infty$. The growth condition  for $f$ (see \eqref{f:growth:2}),  implies that $-\De u=f(x,u)\in L^p(\Om)$ for any $p<r.$ Thus, by the Calderon-Zygmund inequality (see \cite[Theorem 9.14]{G-T}),
$u\in W^{2,p}(\Om),$ for any $p\in(1,r).$ 

\begin{enumerate}
\item[\rm (i)] Assume $r<N.$ Choosing any $p\in(N/2,r)$, by Sobolev embeddings, $u\in W^{1,p^*}(\Om),$ where $\frac1{p^*}:= \frac1{p}- \frac1{N}<\frac1{N}.$ Since $p^*>N,$ $u\in C^{\nu}(\overline{\Om})$ for any $\nu<2-\frac{N}{p}$. Now, from elliptic regularity 
$u\in C^{\nu_0}(\overline{\Om})\cap W^{2,r}(\Om)$ for $\nu_0=2-\frac{N}{r}$.

\item[\rm (ii)] Assume $r=N.$ Choosing any $p\in(N/2,N)$, and reasoning as in (i), $u\in W^{1,p^*}(\Om),$ where $\frac1{p^*}:= \frac1{p}- \frac1{N}<\frac1{N}.$ Also $u\in C^{\nu}(\overline{\Om})$ for any $\nu<1$. Now, from elliptic regularity 
$u\in C^{\nu}(\overline{\Om})\cap W^{2,r}(\Om)$ for  any $\nu<1$.

\item[\rm (iii)] Assume $r>N.$ Choosing any $p\in(N,r)$, and reasoning as above, $u\in C^{1,\nu_0}(\overline{\Om})\cap W^{2,r}(\Om)$ for $\nu_0=1-\frac{N}{r}$.

\item[\rm (iv)] Assume  $r=+\infty.$ Since elliptic regularity and Sobolev embeddings, $u\in C^{1,\nu}(\overline{\Om})\cap W^{2,p}(\Om)$ for any $\nu<1$ and any $p<\infty$.
\end{enumerate}
\end{proof}

\subsection{Estimates of the $L^\infty$-norm of the solutions}
We assume that the non-linearity 
$f$ satisfies the growth condition (H0)', 
and that $\tilde{f}:\R\to (0,+\infty)$  satisfies the following hypothesis:
\begin{enumerate}				
\item[\rm (H1)] there exists a uniform constant $c_0>0$ such that 
\begin{equation}\label{H1}
\limsup_{s\to +\infty}\ \dfrac{\max_{[-s,s]}\, \tilde{f}} {\max\big\{\tilde{f}(-s),\tilde{f}(s)\big\}}\, \le \, c_0.
\footnote{Observe that in particular, if  $\tilde{f}(s)$ is monotone, then (H1) is obviously satisfied with $c_0=1$.} 
\end{equation}	
\end{enumerate}

\bigskip

Under hypothesis  (H0)'-(H1), we establish an estimate for the function $h$ applied to the $L^\infty (\Om)$-norm of any $u\in H_0^{1}(\Om)$ 
solution   to \eqref{eq:ell:pb}, in terms of  their $L^{2^*} (\Om)$-norm. 

From now on, $C$ denotes several constants that may change from line to line, and are independent of $u$. 

Our first main results is the following  theorem.

\begin{thm}
\label{th:apriori:cnys}  
Assume that $f:\Omb \times\R\to \R$  is a Carathéodory function  satisfying  {\rm (H0)'}-{\rm (H1)}.

Then,      for any  $u\in H_0^{1}(\Om)$ weak solution   to \eqref{eq:ell:pb}, the following holds:
\begin{enumerate}				
\item[{\rm (i)}] either there exists a constant $C>0$ such that $\|u\|_{\infty}\le C$, where  $C$  is independent of the  solution $u$,
\item[{\rm (ii)}] either 
for any $\e>0$ there exists a constant $C>0$ such that 
\begin{equation*}%\label{L:inf:h}
h\big(\|u\|_{\infty}\big)
\le C \|a\|_r^{\ A+\e} \ \Big(1+\|u\|_{2^*}\Big)^{\, (2_{N/r}^*-2)(A+\e)},
\end{equation*}
where $h$ is defined by \eqref{def:h}, 
\begin{equation}\label{def:A}
A:= 
\begin{cases}
	1 , &\qq{if}r\le N,\\[.1cm]
	\D1+\frac{2}{N}-\frac{2}{r} ,&\qq{if}r> N,
\end{cases}	
\end{equation}
and $C$ depends only on $\e$, $c_0$ (defined in \eqref{H1}),  $r$, $N$,  and $\Om,$ and it is independent of the  solution $u$.
\end{enumerate}
\end{thm}

As as immediate corollary, we prove that  any sequence of solutions in $H_0^{1}(\Om)$, uniformly bounded in the $ L^{2^*} (\Om)$-norm,  is also uniformly bounded in the $ L^{\infty} (\Om)$-norm.

\begin{cor}\label{cor:apriori:cnys}
Let $f:\Omb \times\R\to \R$ be a   Carathéodory function satisfying {\rm (H0)'}--{\rm (H1)}. 

Let  $\{u_k\}\subset H_0^{1}(\Om)$ be any sequence of solutions to \eqref{eq:ell:pb} such that there exists a constant $C_0>0$ satisfying
\begin{equation*}%\label{L:2*:k:0}
\|u_k\|_{{2^*}} \le C_0.
\end{equation*}

Then, there exists a constant $C>0$ such that
\begin{equation}\label{L:inf:k:0}
\|u_k\|_{\infty} \le C.
\end{equation} 
\end{cor}

\begin{proof}
We reason by contradiction, assuming that \eqref{L:inf:k:0} does not hold. So, at least for a subsequence again denoted as $u_k$, $\|u_k\|_{\infty} \to\infty$  as $k \to \infty.$
Now part (ii) of the Theorem \ref{th:apriori:cnys} implies that 
\begin{equation}\label{h:bdd:uk:0}
h\big(\|u_k\|_{\infty}\big)
\le C.
\end{equation}
From hypothesis (H0)' (see in particular \eqref{h:inf}), for any $\e>0$ there exists $s_1>0$ such that  $h(s)\ge 1/\e$ for any $s\ge s_1,$ and so $h\big(\|u_k\|_{\infty}\big)\ge 1/\e$ for any $k$ big enough. This contradicts \eqref{h:bdd:uk:0},  ending the proof.	
\end{proof}

We next state a straightforward corollary,  assuming that the non-linearity  $\tilde{f}:\R\to (0,+\infty)$  satisfies also the following hypothesis:
\begin{enumerate}				
\item[\rm (H1)'] there exists a uniform constant $c_0>0$ such that 
\begin{equation}\label{H1'}
\sup_{s>0}\ \dfrac{\max_{[-s,s]}\, \tilde{f}} {\max\big\{\tilde{f}(-s),\tilde{f}(s)\big\}}\, \le \, c_0.
\footnote{In particular, if  $\tilde{f}(s)$ is monotone, then (H1)' is  satisfied with $c_0=1$.} 
\end{equation}	
\end{enumerate}

\begin{cor}\label{cor:apriori:cnys:2}
Assume that $f:\Omb \times\R\to \R$  is a Carathéodory function  satisfying  {\rm (H0)'}-{\rm (H1)'}.

Then,      for any  $u\in H_0^{1}(\Om)$ weak solution   to \eqref{eq:ell:pb}, the following holds:
for any $\e>0$ there exists a constant $C>0$ such that 
\begin{equation*}%\label{L:inf:h}
h\big(\|u\|_{\infty}\big)
\le C \|a\|_r^{\ A+\e} \ \Big(1+\|u\|_{2^*}\Big)^{\, (2_{N/r}^*-2)(A+\e)},
\end{equation*}
where $h$ is defined by \eqref{def:h}, $A$ is defined by \eqref{def:A}, 
$C=C(c_0,r,N,\e,|\Om|)$, and $C$ is independent of the  solution $u$.
\end{cor}

\subsection{Proof of Theorem \ref{th:apriori:cnys}}
\label{sec:proof:apriori:cnys}

The arguments of the proof use Gagliardo-Nirenberg interpolation inequality (see \cite{Nirenberg}), and are inspired in the equivalence between uniform $L^{2^*}(\Om)$ {\it a priori} bounds and uniform $L^\infty (\Om)$ {\it a priori} bounds for solutions to subcritical elliptic  equations, see \cite[Theorem 1.2]{Castro-Mavinga-Pardo} for the semilinear case and $f=f(u)$, and \cite[Theorem 1.3]{Mavinga-Pardo_MJM} for the $p$-laplacian  and $f=f(x,u)$.

We first use elliptic regularity and Sobolev embeddings, 
and next, we invoke the Gagliardo-Nirenberg interpolation inequality (see \cite{Nirenberg}).

\begin{proof}[Proof of Theorem \ref{th:apriori:cnys}]
Let   $\{u_k\}\subset H_0^1(\Om)$ be any sequence of  weak solutions to \eqref{eq:ell:pb}. Since Lemma \ref{lem:reg}, in fact $\{u_k\}\subset H_0^1(\Om)\cap L^\infty(\Om)$.

If $\|u_k \|_{\infty}\le C,$ then 
(i)  holds.

Now, we  argue on the  contrary, assuming that there exists a sequence $\|u_k \|_{\infty} \to + \infty$ as $k \to \infty.$

We split the proof in two steps. Firstly, we write an $W^{2,q}$ estimate for $q\in\big(N/2,\min\{r,N\}\big),$ then through Sobolev embeddings we get a $W^{1,q^*}$ estimate with $1/q^*=1/q-1/N<1/N.$ Secondly, we invoke the Gagliardo-Nirenberg interpolation inequality  for the $L^\infty$-norm in terms of its $W^{1,q^*}$-norm and its $L^{2^*}$-norm.

\bigskip

{\it Step 1. $W^{2,q}$ estimates for $q\in\big(N/2,\min\{r,N\}\big)$.}

\medskip

\noindent Let us denote by 
\begin{equation}\label{def:M:k:f}
M_k:=\max\Big\{\tilde{f}\big(-\|u_k\|_{\infty}\big),\tilde{f}\big(\|u_k\|_{\infty}\big)\Big\}
\ge 
(c_0/2)^{-1}\,\max_{[-\|u_k\|_{\infty},\|u_k\|_{\infty}]}\tilde{f},
\end{equation}
where the inequality holds by  hypothesis (H1), see \eqref{H1}.

Let us take $q$ in the interval  $(N/2,N)\cap (N/2,r).$ Growth hypothesis (H0)' (see \eqref{f:growth}), hypothesis (H1) (see \eqref{H1}), and Hölder inequality, yield the following
\begin{align*}%\label{f:q:k:2}
&\displaystyle \int_{\Om}\left|f\big(x,u_k(x)\big)\right|^q\, dx 
\le   \int_{\Om} |a(x)|^q \left(\tilde{f}\big(u_k(x)\big)\right)^{q}\, dx\nonumber\\
&\qquad =   \int_{\Om} |a(x)|^q \left(\tilde{f}\big(u_k(x)\big)\right)^{t} \, 
\left(\tilde{f}\big(u_k(x)\big)\right)^{q-t}\, dx\nonumber\\
&\qquad \le  C
\left[\int_{\Om} |a(x)|^q\ \left(\tilde{f}\big(u_k(x)\big)\right)^{t} \,dx\right]
\ M_k^{\ q-t}
\nonumber\\
&\qquad \le  C   
\left(\int_{\Om} |a(x)|^{qs} \,dx\right)^\frac{1}{s} \left(\int_{\Om} \left(\tilde{f}\big(u_k(x)\big)\right)^{ts'} \,dx\right)^\frac{1}{s'}
\ M_k^{\ q-t}
\nonumber\\
&\qquad \le  C   \|a\|_r^q\ \Big(\|\tilde{f}(u_k) \|_{\frac{2^*}{2_{N/r}^*-1}}\Big)^{t}\
M_k^{\ q-t},
\end{align*}
where $\frac{1}{s}+\frac{1}{s'}=1$, $qs=r$, $C=c_0^{q-t}$  (for $c_0$ defined in \eqref{H1}), and $ts'=\frac{2^*}{2_{N/r}^*-1}$, so 
\begin{align}\label{def:t}
t&:=\frac{2^*}{2_{N/r}^*-1}\left(1-\frac{q}{r}\right)<q \\
&\iff \frac1{q}-\frac1{r}<\frac{2_{N/r}^*-1}{2^*}
=1-\frac1{r}-\frac12+\frac1N \nonumber\\
&\iff \frac1{q}<\frac12+\frac1N\iff q>\frac{2N}{N+2}\ \checkmark \nonumber
\end{align}
since $q>N/2>\frac{2N}{N+2}.$

Now, elliptic regularity and Sobolev embedding imply that
\begin{equation*}%\label{ell:reg}
\|u_k \|_{W^{1,q^*} (\Om)}
\le   C\
\|a\|_r\ \Big(\|\tilde{f}(u_k) \|_{\frac{2^*}{2_{N/r}^*-1}}\Big)^\frac{t}{q}\
M_k^{\ 1-\frac{t}{q}},
\end{equation*}
where $1/q^*=1/q-1/N$, and $C=C(c_0,r,N,q,|\Om|)$ and it is independent of $u.$ Observe that since $q>N/2$, then $q^*>N.$

\bigskip

{\it Step 2. Gagliardo-Nirenberg interpolation inequality.}

\medskip

\noindent Thanks to the Gagliardo-Nirenberg interpolation inequality, there exists a constant $C=C(N,q,|\Om|)$ such that 

\begin{equation*}%\label{Gag-Nir}
\|u_k\|_{\infty}\le C \|\nabla u_k \|_{q^*}^\s \ \|u_k\|_{2^*}^{1-\s}
\end{equation*}
where
\begin{equation}\label{Gag-Nir:2}
\frac{1-\s}{2^*}=\s \left(\frac2{N}-\frac{1}{q}\right).
\end{equation}
Hence
\begin{equation}\label{Gag-Nir:3}
\|u_k\|_{\infty}\le C 
\left[
\|a\|_r\ \Big(\|\tilde{f}(u_k) \|_{\frac{2^*}{2_{N/r}^*-1}}\Big)^{\frac{t}{q}}\
M_k^{\ 1-\frac{t}{q}}	
\right]^\s 	
\ \|u_k\|_{2^*}^{1-\s},
\end{equation}
where $C=C(c_0,r,N,q,|\Om|)$.

From definition of $M_k$ (see \eqref{def:M:k:f}), and definition of $h$ (see \eqref{def:h}),
we deduce that
\begin{equation*}%\label{def:M:k:f:10}
M_k	=\frac{\|u_k\|_{\infty}^{\ 2_{N/r}^*-1}}{h\big(\|u_k\|_{\infty}\big)}.
\end{equation*}
From \eqref{Gag-Nir:2}
\begin{equation}\label{Gag-Nir:4}
\frac{1}{\s}=1+2^*\left(\frac2{N}-\frac{1}{q}\right)=2^*-1-\frac{2^*}{q}
%=2^*\left(\frac12+\frac1N-\frac{1}{q}\right)
=2_{N/q}^*-1.
\end{equation}
Moreover, since definition of  $t$ (see  \eqref{def:t}),
and definition of $2_{N/r}^*$  (see  \eqref{def:2:r}
\begin{align}
\label{def:be:0}
1-\frac{t}{q}&
=\frac{2^*\left(1-\frac{1}{r}\right)-1-2^*\left(\frac1{q}-\frac{1}{r}\right)}{2_{N/r}^*-1}
=\frac{2_{N/q}^*-1}{2_{N/r}^*-1},	
\end{align}
which, joint with \eqref{Gag-Nir:4}, yield
\begin{equation*}\label{def:be:2}
\s\left[1-\frac{t}{q}\right]{(2_{N/r}^*-1)}=1.
\end{equation*}
Now \eqref{Gag-Nir:3} can be rewritten as
\begin{equation*}%\label{Gag-Nir:5}
h\big(\|u_k\|_{\infty}\big)^{\ (1-\frac{t}{q})\s}\le C 
\left[
\|a\|_r\ \Big(\|\tilde{f}(u_k) \|_{\frac{2^*}{2_{N/r}^*-1}}\Big)^{\frac{t}{q}}\		\right]^\s 	
\ \|u_k\|_{2^*}^{1-\s},
\end{equation*}
or equivalently
\begin{equation*}%\label{Gag-Nir:6}
h\big(\|u_k\|_{\infty}\big)\le C 
\|a\|_r^{\ \te} \ \Big(\|\tilde{f}(u_k) \|_{\frac{2^*}{2_{N/r}^*-1}}\Big)^{\te-1}\		
\ \|u_k\|_{2^*}^{\ \vartheta},
\end{equation*}
where 
\begin{align}
\label{def:al}
\te &:=(1-t/q)^{-1}=\frac{2_{N/r}^*-1}{2_{N/q}^*-1},\\
\label{def:al:be}
\vartheta &:=\frac{1-\s}{\s}(1-t/q)^{-1}
=  \te\ (2_{N/q}^*-2),
\end{align}
see \eqref{def:be:0} and \eqref{Gag-Nir:2}.
Observe that since $q< r$, then $\te> 1$. Moreover, since \eqref{def:al:be}, and \eqref{def:al}
\begin{equation}\label{al:be}
\te+\vartheta =\te (2_{N/q}^*-1) =2_{N/r}^*-1.
\end{equation}

Furthermore, from sub-criticallity, see \eqref{f:sub}
$$
\int_\Om |\tilde{f}(u_k) |^{\frac{2^*}{2_{N/r}^*-1}}
\le C \left(1+\int_\Om |u_k|^{{2^*}}\, dx\right),
$$
so
\begin{equation*}%\label{f:u:2*}
\|\tilde{f}(u_k) \|_{\frac{2^*}{2_{N/r}^*-1}}
\le C
\left(1+\|u_k\|_{{2^*}}^{2_{N/r}^*-1}\right).
\end{equation*}
Consequently
\begin{equation*}%\label{Gag-Nir:9}
h\big(\|u_k\|_{\infty}\big)\le C 
\|a\|_r^{\ \te} \ \Big(1+\|u_k\|_{{2^*}}^{\Theta}\Big),
\end{equation*}
with
$$
\Theta:=(2_{N/r}^*-1)(\te-1)+\vartheta=(2_{N/r}^*-2)\te,
$$
where we have used \eqref{al:be}.

Fixed $N>2$ and $r>N/2$, 
the function $q\to \theta=\theta(q)$ for $q\in \big(N/2,\min\{r,N\}\big)$, is decreasing, so 
$$
\inf_{q\in (N/2,\min\{r,N\})} \theta(q)=\theta\big(\min\{r,N\}\big)=A:=\begin{cases}
1,&\text{if}\ r\le N,\\
1+\frac{2}{N}-\frac{2}{r},&\text{if}\ r> N.
\end{cases}
$$
Finally, and since the infimum is not attained in $\big(N/2,\min\{r,N\}\big)$, 
for any $\e>0$,
there exists a constant $C>0$ such that
\begin{equation*}%\label{rad:Gag-Nir:17}
h\big(\|u_k\|_{\infty}\big)\le 
C\ \|a\|_r^{\ A+\e} \ \Big(1+\|u_k\|_{2^*}^{\, (2_{N/r}^*-2)(A+\e)}\Big),
\end{equation*}
where  $A$ is defined by \eqref{def:A}, and $C=C(\e,c_0,r,N,|\Om|),$ ending the proof.
\end{proof}

\section{Radial singular weights}
\label{sec:main:II}
In this section, assuming that $0\in\Om$ and  that   $|f(x,s)|\le |x|^{-\mu}\, \tilde{f}(s)$ for some $\mu\in (0,2)$, we state our second main result concerning weak solutions for singular subcritical non-linearities, see  Theorem \ref{th:rad}. 

\smallskip

First, we also  collect a regularity Lemma for any weak solution to \eqref{eq:ell:pb} with $\tilde{f}(s)$ of polynomial critical growth, according to Caffarelli-Kohn-Nirenberg  inequality.

\begin{lem}[Improved regularity]
\label{lem:reg:rad} 

Assume that $u\in H_0^1(\Om)$ weakly solves \eqref{eq:ell:pb} for a Carathéodory non-linearity $f:\bar{\Om}\times\R\to\R$ with polynomial critical growth
\begin{equation*}%\label{f:growth:4}
|f(x,s)|\le |x|^{-\mu}\,\big(1+|s|^{2_{\mu}^*-1}\big),\qq{with} \mu\in (0,2).
\end{equation*}
Then, the following hold: 
\begin{enumerate}
\item[\rm (i)] If $\mu< 1,$ then $u\in C^{1,\nu}(\overline{\Om})\cap W^{2,p}(\Om)$ for any $p<N/\mu$, and any $\nu<1-\mu$.
\item[\rm (ii)] If $\mu=1,$ then $u\in C^{\nu}(\overline{\Om})\cap W^{2,p}(\Om)$ for  any $p<N$, and $\nu<1$.
\item[\rm (iii)] If $1<\mu<2,$ then $u\in C^{\nu}(\overline{\Om})\cap W^{2,p}(\Om)$ for any $p<N/\mu$, and $\nu<1-\mu$.
\end{enumerate}  

\end{lem} 
\begin{proof}
Let $u\in H_0^1(\Om)$ be a solution to \eqref{eq:ell:pb}. We reason as in Lemma \ref{lem:reg}. 

If $f$ satisfies \eqref{f:sub:rad}, then Caffarelli-Kohn-Nirenberg interpolation inequality (see \cite{Caffarelli-Kohn-Nirenberg}) implies that
assumption \eqref{f:growth:3} is satisfied with
$$
b(x)=\frac{|x|^{-\mu}\,(1+|u|^{2_{\mu}^*-1})}{1+|u|}\le C\, |x|^{-\mu}(1+|u|^{2_{\mu}^*-2})\in L^{N/2}(\Om).
$$
Indeed, since Caffarelli-Kohn-Nirenberg, there exists a constant $C>0$ depending on the parameters $N,$ and $\mu$, such that 
\begin{equation*}%\label{rad:CKN:0}
\big|\,|x|^{\ga}\ u\,\big|_{t}
\le C \|\nabla u\|_{2}^\te \ \|u\|_{2^*}^{1-\te},
\end{equation*}
where $\ga=-\frac{\mu}{2_{\mu}^*-2}=-\frac{\mu(N-2)}{2(2-\mu)}$, $\ t=(2_{\mu}^*-2)\frac{N}{2}=\frac{N(2-\mu)}{N-2}$, and
$\ \frac{1}{t}+ \frac{\ga}{N}=\frac{1}{2^*}= \te \left( \frac12-\frac{1}{N}\right)+(1-\te)\frac{1}{2^*}$, with $\te\in (0,1).$

Consequently, $u\in L^q(\Om)$ for any $q<\infty$. The growth condition  for $f$ (see \eqref{f:growth:rad}-\eqref{f:sub:rad}),  implies that $-\De u=f(x,u)\in L^p(\Om)$ for any $p<N/\mu.$ Thus, by the Calderon-Zygmund inequality (see \cite[Theorem 9.14]{G-T}),
$u\in W^{2,p}(\Om),$ for any $p\in(1,N/\mu ).$ 

\begin{enumerate}
\item[\rm (i)] Assume $\mu< 1.$ Choosing any $p\in(N,N/\mu)$, by elliptic regularity, $u\in W^{2,p}(\Om),$ with $p>N.$ Then $u\in C^{1,\nu}(\overline{\Om})$ for any $\nu<1-\frac{N}{p}$, and finally, 
$u\in C^{1,\nu}(\overline{\Om})\cap W^{2,p}(\Om)$ for any $p<N/\mu$, and any $\nu<1-\mu$.

\item[\rm (ii)] Assume $\mu= 1.$ Choosing any $p\in(N/2,N)$, by elliptic regularity and Sobolev embeddings, $u\in W^{1,p^*}(\Om),$ where $\frac1{p^*}:= \frac1{p}- \frac1{N}<\frac1{N}.$ Also $u\in C^{\nu}(\overline{\Om})$ for any $\nu<1$. Finally 
$u\in C^{\nu}(\overline{\Om})\cap W^{2,p}(\Om)$ for  any $p<N$, and $\nu<1$.

\item[\rm (iii)] Assume $1<\mu<2.$ Choosing any $p\in(N/2,N/\mu)$, and reasoning as above, $u\in W^{1,p^*}(\Om),$ where $\frac1{p^*}:= \frac1{p}- \frac1{N}<\frac1{N}.$ Finally  $u\in C^{\nu}(\overline{\Om})\cap W^{2,p}(\Om)$ for any $p<N/\mu$, and $\nu<1-\mu$.
\end{enumerate}
\end{proof}

\subsection{Estimates of the $L^\infty$-norm of the solutions}

Our second main results is the following  theorem.

\begin{thm}\label{th:rad}
Assume that $f:\Omb \times\R\to \R$  is a Carathéodory function  satisfying  {\rm (H0)} and {\rm (H1)}. 

Then,    for any  $u\in H_0^{1}(\Om)$ solution   to \eqref{eq:ell:pb}, the following holds:
\begin{enumerate}				
\item[{\rm (i)}] either there exists a constant $C>0$ such that $\|u\|_{\infty}\le C$, where  $C$  is independent of the  solution $u$,
\item[{\rm (ii)}] either 
for any $\e>0$ there exists a constant $C>0$ such that 
\begin{align*}%\label{rad:L:inf:h:0}
h\big(\|u\|_{\infty}\big)\le C_\e \ \Big(1+\|u\|_{2^*}\Big)^{\, (2_{\mu}^*-2)(B+\e)}\ , 
\end{align*}
where $h$ is defined by \eqref{rad:def:h},  
\begin{equation}\label{def:B}
B:= 
\begin{cases}
	\D1+\frac{2}{N}-\frac{2\mu}{N} ,&\text{if}\ \mu\in (0,1),\\[.1cm]
	1 ,&\text{if}\ \mu\in [1,2),
\end{cases}	
\end{equation}
and $C$ depends only on $\e$, $c_0$ (defined in \eqref{H1}),  $\mu$, $N$,  and $\Om,$ and it is independent of the  solution $u$.
\end{enumerate}	
\end{thm}

\subsection{Proof of  Theorem \ref{th:rad}}
\label{sec:proof:rad}
\begin{proof}[Proof of Theorem \ref{th:rad}]  
Let   $\{u_k\}\subset H_0^1(\Om)$ be any sequence of  solutions to \eqref{eq:ell:pb}.  Since Lemma \ref{lem:reg:rad}, $\{u_k\}\subset H_0^1(\Om)\cap L^\infty(\Om)$.
If $\|u_k \|_{\infty}\le C,$ then 
(i)  holds.

Now, we  argue on the  contrary, assuming that there exists a sequence $\{u_k\}\subset H_0^1(\Om)$ of  solutions to \eqref{eq:ell:pb}, such that $\|u_k \|_{\infty} \to + \infty$ as $k \to \infty.$ 
%\footnote{
By Morrey's Theorem (see \cite[Theorem 9.12]{Brezis}), observe that also 
\begin{equation}\label{Morrey:th}
\|\na u_k \|_{p} \to + \infty \qq{as} k \to \infty,
\end{equation}
for any $p>N$.
%}

\bigskip

{\it Step 1. $W^{2,q}$ estimates for $q\in\big(N/2,\min\{N,N/\mu\}\big)$.}

\medskip

\noindent As in the proof of Theorem \eqref{th:apriori:cnys}, let us denote by 
\begin{equation}\label{def:M:k:f:rad}
M_k:=\max \Big\{\tilde{f}\big(-\|u_k\|_{\infty}\big),\tilde{f}\big(\|u_k\|_{\infty}\big)\Big\}   \ge 
(c_0/2)^{-1}\max_{[-\|u_k\|_{\infty},\|u_k\|_{\infty}]}\tilde{f},
\end{equation}
where the inequality is due to hypothesis (H1), see \eqref{H1}.

Let us take $q$ in the interval  $(N/2,N)\cap (N/2,N/\mu).$ Using growth hypothesis (H0) (see \eqref{f:growth:rad}), hypothesis (H1) (see \eqref{H1}), and Hölder inequality, we deduce
\begin{align*}%\label{rad:f:q:k:2}
&\displaystyle \int_{\Om}\left|f\big(x,u_k(x)\big)\right|^q\, dx 
\le   \int_{\Om} |x|^{-\mu q} \left(\tilde{f}\big(u_k(x)\big)\right)^{q}\, dx\nonumber\\
&\qquad =   \int_{\Om} |x|^{-\mu q} \left(\tilde{f}\big(u_k(x)\big)\right)^\frac{t}{2_{\mu}^*-1} \, 
\left(\tilde{f}\big(u_k(x)\big)\right)^{q-\frac{t}{2_{\mu}^*-1}}\, dx\nonumber\\
&\qquad \le  C\ 
\left[\int_{\Om} |x|^{-\mu q}\  \big(1+u_k(x)^{t}\big) \,dx\right]
\ M_k^{\ q-\frac{t}{2_{\mu}^*-1}}
\nonumber\\
&\qquad \le C\   \Big(1+\big|\,|x|^{-\ga}\ u_k\,\big|_{t}^{\ t}\Big)\ \
M_k^{\ q-\frac{t}{2_{\mu}^*-1}},
\end{align*}
where $\ga=\frac{\mu q}{t}$,  $t\in \big(0, q\big(2_{\mu}^*-1\big)\big)$, $C=c_0^{q-\frac{t}{2_{\mu}^*-1}}$  (for $c_0$ defined in \eqref{H1}), and where $M_k$ is defined by \eqref{def:M:k:f:rad}.

Combining now elliptic regularity with Sobolev embedding,   we have that
\begin{equation}\label{rad:ell:reg}
\|\nabla u_k \|_{q^*}
\le   C\
\Big(1+\big|\,|x|^{-\ga}\ u_k\,\big|_{t}^{\ t}\Big)^{\,\frac{1}{q}}\
M_k^{\ 1-\frac{t}{q(2_{\mu}^*-1)}},
\end{equation}
where $1/q^*=1/q-1/N$ (since $q>N/2$, then $q^*>N$), and $C=C(N,q,|\Om|).$

\bigskip

{\it Step 2. Caffarelli-Kohn-Nirenberg interpolation inequality.}

\medskip

\noindent Since the Caffarelli-Kohn-Nirenberg interpolation inequality 
for singular weights (see \cite{Caffarelli-Kohn-Nirenberg}), 
there exists a constant $C>0$ depending on the parameters $N,\ q,\ \mu,$ and $t$, such that 
\begin{equation}\label{rad:CKN}
\big|\,|x|^{-\ga}\ u_k\,\big|_{t}
\le C \|\nabla u_k \|_{q^*}^\te \ \|u_k\|_{2^*}^{1-\te},
\end{equation}
where
\begin{equation}\label{rad:CKN:2}
\frac{1}{t}- \frac{\mu q}{N t}= -\te \left( \frac2{N}-\frac{1}{q}\right)+(1-\te)\frac{1}{2^*},\qq{with}\te\in (0,1).
\end{equation}

Substituting now  \eqref{rad:CKN} into \eqref{rad:ell:reg} we can write
\begin{equation*}%\label{rad:CKN:3}
\|\nabla u_k \|_{q^*}\le C
\Big(1+\|\nabla u_k \|_{q^*}^{\te t}\ \|u_k\|_{2^*}^{(1-\te)t}\Big)^\frac{1}{q}\
M_k^{\ 1-\frac{t}{q(2_{\mu}^*-1)}},
\end{equation*}
now, dividing by $\|\nabla u_k \|_{q^*}^{\te t/q}$ and using \eqref{Morrey:th} we obtain
\begin{equation}\label{rad:CKN:4}
\|\nabla u_k \|_{q^*}^{1-\te t/q}\le C
\Big(1+\|u_k\|_{2^*}^{\ \frac{(1-\te)t}{q}}\Big)\
M_k^{\ 1-\frac{t}{q(2_{\mu}^*-1)}}.
\end{equation}

Let us check that
\begin{equation}\label{rad:CKN:4:1}
1-\te\, \frac{t}{q}>0 \qq{for any}t<q\big(2_{\mu}^*-1\big).
\end{equation}
Indeed, observe first that \eqref{rad:CKN:2} is equivalent to 
\begin{equation}\label{rad:CKN:2:1}
\te=\frac{\frac{1}{2^*} -\frac{1}{t}+\frac{\mu q}{N t}}{\frac12+ \frac1{N}-\frac{1}{q}}\ ,
\end{equation}
moreover, from \eqref{rad:CKN:2:1}
\begin{equation}\label{rad:CKN:4:2}
\te\, \frac{t}{q}=\frac{\frac{1}{q}\left(\frac{t}{2^*} -1\right)+ \frac{\mu }{N}}{\frac12+ \frac1{N}-\frac{1}{q}}\ ,
\end{equation}
consequently
\begin{align*}
\te\, \frac{t}{q}<1 &\iff  \frac{1}{q}\left(\frac{t}{2^*} -1\right)+ \frac{\mu }{N}< \frac12+ \frac1{N}-\frac{1}{q}\\
&\iff
\frac{1}{q}\frac{t}{2^*}< 
\frac12+ \frac1{N}-\frac{\mu }{N}\pm 1\\
&\iff
\frac{t}{q}<2^*\left(1- \frac{\mu }{N }\right)-2^*\left(\frac12- \frac1{N}\right)=2_{\mu}^*-1
\\
&\iff
t<q\big(2_{\mu}^*-1\big),
\end{align*}
so, \eqref{rad:CKN:4:1} holds.

Consequently,
\begin{equation}\label{rad:CKN:5}
\|\nabla u_k \|_{q^*}\le C
\Big(1+ \|u_k\|_{2^*}^{\ \frac{(1-\te)t}{q-\te t}}\Big)\
M_k^{\ \big(1-\frac{t}{q(2_{\mu}^*-1)}\big)(1-\te t/q)^{-1}}.
\end{equation}

\bigskip

{\it Step 3. Gagliardo-Nirenberg interpolation inequality.}

\medskip

\noindent Thanks to the Gagliardo-Nirenberg interpolation inequality 
(see \cite{Nirenberg}), there exists a constant $C=C(N,q,|\Om|)$ such that 

\begin{equation}\label{rad:Gag-Nir}
\|u_k\|_{\infty}
\le C \|\nabla u_k \|_{q^*}^\s \ \|u_k\|_{2^*}^{1-\s},
\end{equation}
where
\begin{equation}\label{rad:Gag-Nir:2}
\frac{1-\s}{2^*}=\s \left(\frac2{N}-\frac{1}{q}\right).
\end{equation}
Hence, substituting \eqref{rad:CKN:5} into \eqref{rad:Gag-Nir} we deduce
\begin{equation}\label{rad:Gag-Nir:12}
\|u_k\|_{\infty}\le C 
\ \Big(1+ \|u_k\|_{2^*}^{\ \s\,\frac{(1-\te)t}{q-\te t}
+1-\s}\Big)
\ M_k^{\ \s \big(1-\frac{t}{q(2_{\mu}^*-1)}\big)(1-\te t/q)^{-1}}		
.
\end{equation}

\medskip

From definition of $M_k$ (see \eqref{def:M:k:f}), and 
of $h$  (see \eqref{rad:def:h}),
we obtain
\begin{equation}\label{rad:def:M:k:f:10}
M_k	=\frac{\|u_k\|_{\infty}^{2_{\mu}^*-1}}{h\big(\|u_k\|_{\infty}\big)}.
\end{equation}
From \eqref{rad:Gag-Nir:2}
\begin{equation}\label{rad:Gag-Nir:13}
\frac{1}{\s}=1+2^*\left(\frac2{N}-\frac{1}{q}\right)
=2_{N/q}^*-1.
\end{equation}
From \eqref{rad:CKN:4:2}, 
we deduce
\begin{align}\label{rad:CKN:4:3}
1-\te\, \frac{t}{q}&=\frac{\frac12+ \frac1{N}-\frac{t}{2^*q} - \frac{\mu}{N}\pm 1}{\frac12+ \frac1{N}-\frac{1}{q}}\\
&=\frac{\left(1- \frac{\mu}{N}\right)-\frac1{2^*} -\frac{t}{2^*q} } {\frac12+ \frac1{N}-\frac{1}{q}}
=\frac{2_{\mu}^*-1-\frac{t}{q}}{2_{N/q}^*-1}
\nonumber
\ ,
\end{align}
where we have used that, by definition of $2_{\mu}^*$ (see \eqref{def:2*:mu}), $\frac{2_{\mu}^*}{2^*}=1- \frac{\mu}{N}$. 

Moreover, since  \eqref{rad:CKN:4:3},
\begin{align}\label{rad:def:be}
&\left(1-\frac{t}{q(2_{\mu}^*-1)}\right){\big(2_{\mu}^*-1\big)}
\frac1{(1-\te t/q)}
%\nonumber
\\ 
&\qquad =\left( 2_{\mu}^*-1-\frac{t}{q}\right)
\frac1{(1-\te t/q)}=2_{N/q}^*-1.
\nonumber	
\end{align}
Taking into account \eqref{rad:Gag-Nir:13} and  \eqref{rad:def:be} we obtain
\begin{equation}\label{1}
\s\Big(1-\frac{t}{q(2_{\mu}^*-1)}\Big){\big(2_{\mu}^*-1\big)}(1-\te t/q)^{-1}=1.
\end{equation} 
Consequently, since \eqref{rad:def:M:k:f:10}, and \eqref{1}, we can rewrite \eqref{rad:Gag-Nir:12} in the following way
\begin{equation}\label{rad:Gag-Nir:14}
h\big(\|u_k\|_{\infty}\big)^{\ \frac1{2_{\mu}^*-1}}\le C \ \Big(1+ 
\|u_k\|_{2^*}^{\ \s\,\frac{(1-\te)t/q}{1-\te t/q}+1-\s}\Big),
\end{equation}
or equivalently
\begin{equation}\label{rad:Gag-Nir:15}
h\big(\|u_k\|_{\infty}\big)\le C\ \Big(1+  \|u_k\|_{2^*}^{\ \Theta}\Big),
\end{equation}
where
\begin{align*}
\Theta &:=
\big(2_{\mu}^*-1\big)\left[1+\s\,\frac{t/q-1}{1-\te t/q}\right].	
\end{align*}

\medskip

Since \eqref{1}, $\s(1-\te t/q)^{-1}=(2_{\mu}^*-1 -\frac{t}{q})^{-1}$, and substituting it into the above equation we obtain
\begin{align*}
\Theta &=\big(2_{\mu}^*-1\big)\left( \frac{2_{\mu}^*-2}{2_{\mu}^*-1 -\frac{t}{q}}\right).
\end{align*}
Fixed $N>2$ and $\mu\in (0,2)$, the function $(t,q)\to \Theta=\Theta(t,q)$ for $(t,q)\in \big(0, q(2_{\mu}^*-1)\big)\times \big(N/2,\min\{N,N/\mu\}\big)$, is increasing in $t$ and decreasing in $q$.

For $\mu\in[1,2)$,  $\min\{N,N/\mu\}=N/\mu$. If $q_k\to N/\mu$, equation \eqref{rad:CKN:2} with $q=q_k$,  $\te=\te_k<1$ and an arbitrary $t\in \big(0, (2_{\mu}^*-1)N/\mu\big)$ fixed,   yields $\te_k\to\frac{1}{2_{\mu}^*-1}<1$ (since $\mu<2$). 
Hence, when $\mu\in[1,2)$,
$$
\underset{t\in\big(0, (2_{\mu}^*-1)\frac{N}{\mu}\big),\, q\in 
\big(\frac{N}{2},\frac{N}{\mu}\big)}{\inf}\ \Theta(t,q)
=\Theta\bigg(0,\frac{N}{\mu}\bigg)=2_{\mu}^*-2.
$$

On the other hand, for $\mu\in(0,1)$,  $\min\{N,N/\mu\}=N$. If $q_k\to N,$ equation \eqref{rad:CKN:2} with $q=q_k$, $\te=\te_k>0$ and $t$ fixed,  yields $\te_k\to\frac{2}{2^*}-\frac{2(1-\mu)}{t}\ge 0$, so $t\ge 2^*(1-\mu)$. Hence, when $\mu\in(0,1)$,
$$
\underset{t\in[2^*(1-\mu), (2_{\mu}^*-1)N),\, 
q\in \big(\frac{N}{2},N\big)}{\inf}\,  \Theta(t,q)
=\Theta(2^*(1-\mu),N)=(2_{\mu}^*-2)B,
$$
where $B$ is defined by \eqref{def:B}.

Since the infimum is not attained, for any $\e>0$,
there exists a constant $C=C(\e,c_0,\mu,N,\Om)$ such that 
\begin{equation}\label{rad:Gag-Nir:16}
h\big(\|u_k\|_{\infty}\big)\le C \ 
\Big(1+\|u_k\|_{2^*}^{\, (2_{\mu}^*-2)(B+\e)}\Big) ,
\end{equation}
which ends the proof.
\end{proof}

\section{Acknowledgments}
The author is supported by Grant PID2019-103860GB-I00,  MICINN,  Spain and Grupo de Investigación CADEDIF 920894, UCM.

\def\cprime{$'$}


\begin{thebibliography}{99}

\bibitem{Bidaut-Veron}
M.-F. Bidaut-Veron.
\newblock Local behaviour of the solutions of a class of nonlinear elliptic
systems.
\newblock {\em Adv. Differential Equations}, 5(1-3):147--192, 2000.

\bibitem{Brezis}
H.~Brezis.
\newblock {\em Functional analysis, {S}obolev spaces and partial differential
equations}.
\newblock Universitext. Springer, New York, 2011.

\bibitem{Brezis-Kato}
H.~Br{\'e}zis and T.~Kato.
\newblock Remarks on the {S}chr\"odinger operator with singular complex
potentials.
\newblock {\em J. Math. Pures Appl. (9)}, 58(2):137--151, 1979.

\bibitem{Caffarelli-Kohn-Nirenberg}
L.~Caffarelli, R.~Kohn, and L.~Nirenberg.
\newblock First order interpolation inequalities with weights.
\newblock {\em Compositio Math.}, 53(3):259--275, 1984.

\bibitem{Caldiroli-Malchiodi}
P.~Caldiroli and A.~Malchiodi.
\newblock Singular elliptic problems with critical growth.
\newblock {\em Comm. Partial Differential Equations}, 27(5-6):847--876, 2002.

\bibitem{Castro-Mavinga-Pardo}
A.~Castro, N.~Mavinga, and R.~Pardo.
\newblock Equivalence between uniform {$L^{2^\ast}(\Omega)$} a-priori bounds
and uniform {$L^\infty(\Omega)$} a-priori bounds for subcritical elliptic
equations.
\newblock {\em Topol. Methods Nonlinear Anal.}, 53(1):43--56, 2019.

\bibitem{Castro_Pardo_RMC}
A.~Castro and R.~Pardo.
\newblock A priori bounds for positive solutions of subcritical elliptic
equations.
\newblock {\em Rev. Mat. Complut.}, 28(3):715--731, 2015.

\bibitem{Castro-Pardo_DCDS}
A.~Castro and R.~Pardo.
\newblock A priori estimates for positive solutions to subcritical elliptic
problems in a class of non-convex regions.
\newblock {\em Discrete Contin. Dyn. Syst. Ser. B}, 22(3):783--790, 2017.

\bibitem{Clapp-Pardo-Pistoia-Saldana}
M.~Clapp, R.~Pardo, A.~Pistoia, and A.~Salda\~{n}a.
\newblock A solution to a slightly subcritical elliptic problem with non-power
nonlinearity.
\newblock {\em J. Differential Equations}, 275:418--446, 2021.

\bibitem{Damascelli-Pardo}
L.~Damascelli and R.~Pardo.
\newblock A priori estimates for some elliptic equations involving the
{$p$}-{L}aplacian.
\newblock {\em Nonlinear Anal. Real World Appl.}, 41:475--496, 2018.

\bibitem{deFigueiredo-Lions-Nussbaum}
D.~G. de~Figueiredo, P.-L. Lions, and R.~D. Nussbaum.
\newblock A priori estimates and existence of positive solutions of semilinear
elliptic equations.
\newblock {\em J. Math. Pures Appl. (9)}, 61(1):41--63, 1982.

\bibitem{Dupaigne-Ponce}
L.~Dupaigne and A.~C. Ponce.
\newblock Singularities of positive supersolutions in elliptic {PDE}s.
\newblock {\em Selecta Math. (N.S.)}, 10(3):341--358, 2004.

\bibitem{Gidas-Spruck2}
B.~Gidas and J.~Spruck.
\newblock Global and local behavior of positive solutions of nonlinear elliptic
equations.
\newblock {\em Comm. Pure Appl. Math.}, 34(4):525--598, 1981.

\bibitem{Gidas-Spruck}
B.~Gidas and J.~Spruck.
\newblock A priori bounds for positive solutions of nonlinear elliptic
equations.
\newblock {\em Comm. Partial Differential Equations}, 6(8):883--901, 1981.

\bibitem{G-T}
D.~Gilbarg and N.~S. Trudinger.
\newblock {\em Elliptic partial differential equations of second order}, volume
224 of {\em Grundlehren der Mathematischen Wissenschaften [Fundamental
Principles of Mathematical Sciences]}.
\newblock Springer-Verlag, Berlin, second edition, 1983.

\bibitem{Jannelli-Solimini}
E.~Jannelli and S.~Solimini.
\newblock Critical behaviour of some elliptic equations with singular
potentials.
\newblock {\em Rapporto no. 41/96, Dipartimento Di Mathematica, Universita
degli Studi di Bari, 70125 Bari, Italia.}

\bibitem{Joseph-Lundgren}
D.~D. Joseph and T.~S. Lundgren.
\newblock Quasilinear {D}irichlet problems driven by positive sources.
\newblock {\em Arch. Rational Mech. Anal.}, 49:241--269, 1972/73.

\bibitem{Leray-Schauder}
J.~Leray and J.~Schauder.
\newblock Topologie et \'{e}quations fonctionnelles.
\newblock {\em Ann. Sci. \'{E}cole Norm. Sup. (3)}, 51:45--78, 1934.

\bibitem{Mavinga-Pardo17_JMAA}
N.~Mavinga and R.~Pardo.
\newblock A priori bounds and existence of positive solutions for semilinear
elliptic systems.
\newblock {\em J. Math. Anal. Appl.}, 449(2):1172--1188, 2017.

\bibitem{Mavinga-Pardo_MJM}
N.~Mavinga and R.~Pardo.
\newblock Equivalence between uniform {$L^{p^*}$} a priori bounds and uniform
{$L^\infty$} a priori bounds for subcritical {$p$}-{L}aplacian equations.
\newblock {\em Mediterr. J. Math.}, 18(1):Paper No. 13, 24, 2021.

\bibitem{Moser}
J.~Moser.
\newblock A new proof of {D}e {G}iorgi's theorem concerning the regularity
problem for elliptic differential equations.
\newblock {\em Comm. Pure Appl. Math.}, 13:457--468, 1960.

\bibitem{Nirenberg}
L.~Nirenberg.
\newblock On elliptic partial differential equations.
\newblock {\em Ann. Scuola Norm. Sup. Pisa Cl. Sci. (3)}, 13:115--162, 1959.

\bibitem{Pardo19a}
R.~Pardo.
\newblock On the existence of a priori bounds for positive solutions of
elliptic problems, {I}.
\newblock {\em Rev. Integr. Temas Mat.}, 37(1):77--111, 2019.

\bibitem{Pardo19b}
R.~Pardo.
\newblock On the existence of a priori bounds for positive solutions of
elliptic problems, {II}.
\newblock {\em Rev. Integr. Temas Mat.}, 37(1):113--148, 2019.

\bibitem{Pardo-Sanjuan}
R.~Pardo and A.~Sanju\'{a}n.
\newblock Asymptotic behavior of positive radial solutions to elliptic
equations approaching critical growth.
\newblock {\em Electron. J. Differential Equations}, pages Paper No. 114, 17,
2020.

\bibitem{Struwe}
M.~Struwe.
\newblock {\em Variational methods}, volume~34 of {\em Ergebnisse der
Mathematik und ihrer Grenzgebiete. 3. Folge. A Series of Modern Surveys in
Mathematics [Results in Mathematics and Related Areas. 3rd Series. A Series
of Modern Surveys in Mathematics]}.
\newblock Springer-Verlag, Berlin, fourth edition, 2008.
\newblock Applications to nonlinear partial differential equations and
Hamiltonian systems.

\end{thebibliography}
\end{document}